\documentclass{amsart}
\usepackage{amssymb}
\usepackage{amsmath}
\usepackage{textcomp}
\makeindex

\newcommand{\ba}{\begin{array}}
\newcommand{\eea}{\end{eqnarray}}
\newcommand{\ea}{\end{array}}

\newtheorem{definition}{Definition}[section]
\newtheorem{theorem}[definition]{Theorem}

\newtheorem{example}[definition]{Example}

\usepackage[matrix,arrow,curve]{xy}
\usepackage[utf8x]{inputenc}
\usepackage{tikz}

\begin{document}
\title[Seifert's Conjecture For Almost Symplectic Foliations]{Seifert's Conjecture For Almost Symplectic Foliations}
\author[S. Mukherjee]{Sauvik Mukherjee}


\begin{abstract} We disproving Seifert's conjecture for almost symplectic foliations with co-dimension bigger or equal to $3$.\\
% we show that the closedness of the leaves is not an obstruction for the $h$-principle (the existence part only that too without homotopy) for symplectic foliations on closed manifolds. 
\end{abstract}
\maketitle

\section{introduction} Let $Gr_{2n}(M)\stackrel{\pi}{\to} M$ be the Grassmann bundle on $M^{2n+q}$, i.e $\pi^{-1}(x)=Gr_{2n}(T_xM)$. Identify $Fol_q(M)$, the space of codimension-$q$ foliations on $M$ as a subspace of $\Gamma(Gr_{2n}(M))$, where the section space $\Gamma(Gr_{2n}(M))$ is given the $C^{\infty}$ topology. As usual $\Omega^k(M)$ be the space of $k$-forms on $M$ with $C^{\infty}$-topology.   Define \[\Delta_q(M)\subset Fol_q(M)\times \Omega^2(M)\] \[\Delta_q(M)=\{(\mathcal{F},\omega): \omega^n_{\mid \mathcal{F}}\neq 0\}\]

  Now we state the main result of this paper. Let $\bar{\Delta}_q(M)$ be the space of all pairs $(\mathcal{F},\omega)\in \Delta_q(M)$ such that $\mathcal{F}$ does not have a closed leaf.
 
 \begin{theorem}(Main Theorem)
  \label{Main}
  Let $q\geq 3,\ n\geq 2$ and $(\mathcal{F}_0,\omega_0)$ be a pair in $\Delta_q(M)$ then there exists a pair $(\mathcal{G},\gamma)\in \bar{\Delta}_q(M)$.  
 \end{theorem}
 
 In the absence of a $2$-form and for foliations with co-dimension bigger or equal to $3$ it has been studied by Schweitzer in \cite{SchWeitzer}. He studied this in the light of Seifert's conjecture which states the following\\
 
 {\bf Seifert's Conjecture:} Every non-vanishing vector field on the three sphere $\mathbb{S}^3$ has a closed integral curve.\\
 
 He considers the general case namely if smooth foliation on a closed manifold admits a closed leaf. He disproves this result in case of foliations with co-dimension bigger or equal to three. \\
 
  We end this section with the definition of a minimal set which is a key ingredient in our proof. \\
 
 \begin{definition}(\cite{SchWeitzer})
  \label{Minimal set}
  Let $M$ be a manifold with a foliation $\mathcal{F}$, a set $S\subset M$ is called saturated with respect to $\mathcal{F}$ if for all $x\in S$ the leaf through $x$ is contained in $S$. $S$ is called minimal if it is nonempty, closed, saturated and minimal in a sense that it does not contain any proper subset with all these properties. 
 \end{definition}
\begin{example}
 Consider the torus $\mathbb{T}^2$ with coordinates $(a,b)$. Let $\mathcal{F}$ be the foliation on $\mathbb{T}^2$ defined by the vector field $\partial_a+c\partial_b$, where $c$ is irrational. Then $\mathbb{T}^2$ is the minimal set for $\mathcal{F}$.
\end{example}

\section{passing to the local model}\label{Passing to the local model} In our proof we shall follow the methods in \cite{SchWeitzer}. In \cite{Wilson} the following result has been proved and which allows us to reduce the original problem of proving \ref{Main} to splitting the leaves of the product foliation on the products of discs.

\begin{theorem}(\cite{Wilson})
 \label{Wilson}
 Let $M^{2n+q}$ be a smooth manifold together with a codimension-$q$ foliation $\mathcal{F}$ on it then there exists a family of embeddings $f_{\lambda}:\mathbb{D}^q \times \mathbb{D}^{2n}\to M,\ \lambda \in \Lambda$ such that the family of sets $\{f_{\lambda}(\mathbb{D}^q \times \mathbb{D}^{2n} ):\lambda \in \Lambda\}$ is locally finite, mutually disjoint. Moreover $f_{\lambda}^{-1}(\mathcal{F})$ has leaves $\{x\}\times \mathbb{D}^{2n}$ and each leaf of $\mathcal{F}$ intersects atleast one of $f_{\lambda}(\mathbb{D}_{1/2}^q \times \mathbb{D}^{2n})$, where $\mathbb{D}_{1/2}^q$ is the $q$-disc of radious $1/2$. 
 \end{theorem}
We shall consider $(f_{\lambda}^{-1}(\mathcal{F}),f_{\lambda}^*\omega_0)$ on $\mathbb{D}^q \times \mathbb{D}^{2n}$. Obviously $(f_{\lambda}^{-1}(\mathcal{F}),f_{\lambda}^*\omega_0)\in \Delta_q(\mathbb{D}^q \times \mathbb{D}^{2n})$.\\

% As in \cite{SchWeitzer} we shall split the leaves of $f_{\lambda}^{-1}(\mathcal{F})$ whose leaves are $\{x\}\times \mathbb{D}^{2n},\ x\in \mathbb{D}^q$ by a concordance. This concordance will be fixed near the boundary of $\mathbb{D}^q \times \mathbb{D}^{2n}$. Then we shall construct a homotopy of two forms $\omega_t$ starting from $f_{\lambda}^{-1}(\omega_0)=\omega'$ satisfying the conditions of \ref{Main} and $\omega_t=\omega',\ on\ Op(\partial(\mathbb{D}^q \times \mathbb{D}^{2n}))$. This would suffice according to \ref{Wilson}.

\section{$h$-principle for the 2-form} According to section-\ref{Passing to the local model} it is enough to consider $(\mathcal{F}_0,\omega_0)$ on $\mathbb{D}^q \times \mathbb{D}^{2n}$ where $\mathcal{F}_0$ is given by  $\{x\}\times \mathbb{D}^{2n},\ x\in \mathbb{D}^q$ and $\omega_0$ is such that $(\omega_0)^n_{\mid T\mathcal{F}_0}\neq 0$. In the first step we change the $2$-form $\omega_0$ to $\omega_1$ such that 
\begin{enumerate}
\item $(\mathcal{F}_0,\omega_1)\in \Delta_q(M)$\\
\item outside the sets $\{f_{\lambda}(\mathbb{D}^q \times \mathbb{D}^{2n} ):\lambda \in \Lambda\}$ (\ref{Wilson}) we have $\omega_1 =\omega_0$\\
\item  the restriction of $\omega_1$ to $\mathbb{D}^q_{\varepsilon}\times \mathbb{D}_{\varepsilon}^{2n}$ is $0 \oplus (\Sigma_1^n\pm (dx_i\wedge dy_i))$\\
\end{enumerate}

Let us prove this fact. Let $(x_i,y_i):i=1,...,n$ be the coordinates in $\mathbb{D}^{2n}$ and the coordinates in $\mathbb{D}^{q}$ be $(z_j):j=1,...,q$. So any two form, in particular $\omega_0$ is of the form \[\omega_0=\Sigma f_{i,j} dx_i\wedge dy_j+\Sigma_{i<j} g_{i,j} dx_i\wedge dx_j+\Sigma_{i<j} h_{i,j} dy_i\wedge dy_j+ \Omega\] where $\Omega$ is a two form of the form $\Omega=\Sigma_{i,j} F_{i,j} dx_i\wedge dz_j+ \Sigma_{i,j} G_{i,j} dy_i\wedge dz_j$.\\

So we get $(\omega_0^n)_{\mid{T\mathcal{F}_0}}$ is of the form $(A_1+...+A_k)dx_1\wedge dy_1\wedge dx_2 \wedge dy_2\wedge ...\wedge dx_n\wedge dy_n$ where each $A_i$ is a product of $f_{i,j},\ g_{i,j}\ and\ h_{i,j}$. Set $f_i=f_{i,i}$. One of the $A_i$'s say $A_1$ is of the form $f_1...f_n$.\\

As $(\omega_0)^n_{\mid T\mathcal{F}_0}\neq 0$ we may assume with out loss of generality that $(\omega_0)^n_{\mid T\mathcal{F}_0}$ evaluated at $0\in \mathbb{D}^{2n}\times \mathbb{D}^q$ is positive. So $(A_1+...+A_k)(0)>0$ and hence $(A_1+...+A_k)$ remains positive in a neighborhood of $0\in \mathbb{D}^{2n}\times \mathbb{D}^q$. One of $A_i(0)$'s is the biggest positive among all $A_i(0)$'s. Say this be $A_1(0)$ otherwise we may get some combination other than $0 \oplus (\Sigma \pm dx_i\wedge dy_i)$ say \[0\oplus (dx_1\wedge dx_2+dy_1\wedge dy_2+\Sigma_3^n(dx_i\wedge dy_i))\] in place of $\omega_1$. But this does not make any difference as we can work with $0\oplus (dx_1\wedge dx_2+dy_1\wedge dy_2+\Sigma_3^n(dx_i\wedge dy_i))$ in place of $\omega_1$ for the rest of the argument.\\

Now as $A_1(0)>0$ is the biggest positive number among all $A_i(0)$, so it remains so in a neighborhood of $0\in \mathbb{D}^{2n}\times \mathbb{D}^q$ say $\mathbb{D}^{2n}_{\varepsilon_1}\times \mathbb{D}^q_{\varepsilon_1}$ for some $\varepsilon_1>0$. Let \[B:\mathbb{D}^{2n}\times \mathbb{D}^q\to \mathbb{R}\] be a smooth function such that $B=0$ on $\mathbb{D}^{2n}_{\varepsilon_1/3}\times \mathbb{D}^q_{\varepsilon_1/3}$ and $B=1$ outside $\mathbb{D}^{2n}_{\varepsilon_1}\times \mathbb{D}^q_{\varepsilon_1}$.\\

Now if we replace $\omega_0$ by \[\bar{\omega}_0=\Sigma f_i dx_i\wedge dy_j+\Sigma_{i<j}B g_{i,j} dx_i\wedge dx_j+\Sigma_{i<j}B h_{i,j} dy_i\wedge dy_j+B \Omega\] Then $(\bar{\omega}_0)^n_{\mid T\mathcal{F}_0}$ is equal to $(\omega_0)^n_{\mid T\mathcal{F}_0}$ outside $\mathbb{D}^{2n}_{\varepsilon_1}\times \mathbb{D}^q_{\varepsilon_1}$ and on $\mathbb{D}^{2n}_{\varepsilon_1/3}\times \mathbb{D}^q_{\varepsilon_1/3}$ it is equal to $(f_1...f_n)dx_1\wedge dy_1\wedge dx_n\wedge dy_n$. As $(f_1...f_n)>0$ on $\mathbb{D}^{2n}_{\varepsilon_1/3}\times \mathbb{D}^q_{\varepsilon_1/3}$ so we can make $f_i$'s either $+1$ or $-1$ on $\mathbb{D}^{2n}_{\varepsilon}\times \mathbb{D}^q_{\varepsilon}$ for $0<\varepsilon< \varepsilon_1/3$.\\

%{\bf Construction:} Observe that the restriction of $\omega_0$ to $\mathbb{D}^q \times \mathbb{D}_{\varepsilon}^{2n}$ is of the form \[(\Sigma_1^nf_idx_i\wedge dy_i+\Sigma_{i,j}g_{ij}dx_i\wedge dx_{j}+ \Sigma_{i,j}h_{ij}dy_i\wedge dy_j)+\Omega\] where $f_i$'s are always non-zero and is either positive or negative through out $\mathbb{D}_{\varepsilon}^{2n}\times \mathbb{D}^q$. So we can homotop $\omega_0$ through leafwise non-degenerate $2$-forms to the desired form by making $f_i$'s $+1$ or $-1$ depending on the sign of $f_i$'s and making $g_{ij}$'s,$h_{ij}$'s and $\Omega$ zero. \\
  
 % So without loss of generality we can keep our model $\mathbb{D}^q\times \mathbb{D}^{2n}$ instead of the modified one $\mathbb{D}^q\times \mathbb{D}_{\varepsilon}^{2n}$.\\

  \section{final step}  As $q\geq 3$ there exists an embedding $\mathbb{T}^2\hookrightarrow \mathbb{D}^q$. Let $N(\mathbb{T}^2)$ be a small tubular neighborhood of $\mathbb{T}^2$ in $\mathbb{D}^q$. So there exists a diffeomorphism $e:\mathbb{T}^2\times \mathbb{D}^{q-2}\to N(\mathbb{T}^2)$. So now consider $(e\times id_{\mathbb{D}^{2n}})^*\omega_1=0\oplus (\Sigma_1^n\pm dx_i\wedge dy_i)$. Call this new form $\omega'_1$ on $\mathbb{T}^2\times \mathbb{D}^{q-2}\times \mathbb{D}^{2n}$. Also set $\mathcal{F}'_1=(e\times id_{\mathbb{D}^{2n}})^{-1}\mathcal{F}_0$. So obviously $(\mathcal{F}'_1,\omega'_1)\in \Delta_q(\mathbb{T}^2\times \mathbb{D}^{q-2}\times \mathbb{D}^{2n})$ and hence this is the model now.\\

Now let $Z$ be a vector field on $\mathbb{T}^2$ for which $\mathbb{T}^2$ is the minimal set (\ref{Minimal set}). Let $Z'=(Z,0)$ be the corresponding vector field on $\mathbb{T}^2\times \mathbb{D}^{q-2}$. Let $\psi:\mathbb{T}^2\times \mathbb{D}^{q-2} \times [0,1]\to [0,1]$ be a smooth function with compact support and supported in the interior of $\mathbb{T}^2\times \mathbb{D}^{q-2}\times [0,1]$ such that $\psi^{-1}(1)=\mathbb{T}^2\times \mathbb{D}_{1/2}^{q-2} \times \{1/2\}$. Set \[X_1=(1-\psi).\partial_s+\psi.Z'\] where $s$ is the variable on $[0,1]$. Observe that as $\psi$ is compactly supported in the interior of $\mathbb{T}^2\times \mathbb{D}^{q-2}\times [0,1]$, so $X_1=\partial_s$ near $\mathbb{T}^2\times \mathbb{D}^{q-2}\times \{1\}$. So we can extend $X_1$ to all of $\mathbb{T}^2\times \mathbb{D}^{q-2} \times [0,\infty)$. We shall denote this new extended vector field by $X_1$ itself.\\

Consider the map $\rho:\mathbb{D}^{2n}\to [0,\infty)$ given by $\rho(z)=\frac{1}{|z|}-1$. So \[id\times \rho: \mathbb{T}^2\times \mathbb{D}^{q-2}\times \mathbb{D}^{2n}\to \mathbb{T}^2\times \mathbb{D}^{q-2}\times [0,\infty)\] is a submersion on $\mathbb{T}^2\times \mathbb{D}^{q-2}\times (\mathbb{D}^{2n}-\{0\})$. Observe that $d\rho_z=-\frac{1}{|z|^3}(x_1,y_1,...,x_n,y_n)$, where $0\neq z=(x_1,y_1,...,x_n,y_n)$. Let $\tilde{\mathcal{F}}_1$ be the foliation on $\mathbb{T}^2\times \mathbb{D}^{q-2}\times [0,\infty)$ defined by the vector field $X_1$. Now $(id\times \rho_0)^{-1}\tilde{\mathcal{F}}_1$ induces a foliation on $\mathbb{T}^2\times \mathbb{D}^{q-2}\times (\mathbb{D}^{2n}-\{0\})$ which we denote by $\mathcal{F}''_1$ , where $\rho_0$ is the restriction of $\rho$ to $ \mathbb{D}^{2n}-\{0\}$. \\

%This foliation agrees with $\mathcal{F}'_1$ near $\mathbb{T}^2\times \mathbb{D}^{q-2}\times \{0\}$ and hence extends to a foliation $\mathcal{F}'_2$ on $\mathbb{T}^2\times \mathbb{D}^{q-2}\times \mathbb{D}^{2n}$. The homotopy from $\mathcal{F}'_1$ to $\mathcal{F}'_2$ is defined by pulling back the codimension-$(q+1)$ foliation defined by $X_t$ on $\mathbb{T}^2\times \mathbb{D}^{q-2}\times [0,\infty)\times [0,1]$.\\

Now we shall construct the two form $\omega''_1$ from $\omega'_1$ and satisfying the conditions of \ref{Main} in the complement of the origin $\{z=0\}$. \\   

Observe that $Y=\frac{1}{2}\Sigma_1^n(x_i\partial_{x_i}+y_i\partial_{y_i})$ is a Liouville vector field of $\Sigma_1^n(\pm dx_i\wedge dy_i)$ on $\mathbb{D}^{2n}$, i.e, $d(i_{Y}(\Sigma_1^n(\pm dx_i\wedge dy_i)))=\Sigma_1^n(\pm dx_i\wedge dy_i)$, where $i_Y$ is the contraction by the vector field $Y$. Set $\alpha=i_{Y}(\Sigma_1^n(\pm dx_i\wedge dy_i))$. Let $H$ be any hyperplane in $\mathbb{R}^{2n}$ transversal to $Y$ then by $1.4.5$ of \cite{Geiges} we have $(\alpha \wedge (d\alpha)^n)_{\mid H}\neq 0$. \\

Observe $(d(id\times \rho))(T\mathcal{F}''_1)=X_1$. Let $A=(a_1,b_1,...,a_n,b_n)\in \mathbb{R}^{2n}$ then for $w\in \mathbb{T}^2\times \mathbb{D}^{q-2}$ and $z=(x_1,y_1,...,x_n,y_n)\in \mathbb{D}^{2n}$, \[(T\mathcal{F}''_1)_{(w,z)}=\{A+\psi Z':-\frac{1}{|z|^3}\Sigma(a_ix_i+b_iy_i)=(1-\psi)\}\] To see this let $A+B\in T_w(\mathbb{T}^2\times \mathbb{D}^2)\times T_z\mathbb{D}^{2n}$ and \[d(id\times \rho)_{(w,z)}(A+B)=d\rho_z(A)+B=X_1\]
So we get $B=\psi Z'$ and $d\rho_z(A)=(1-\psi)$.\\

 Now $T\mathcal{F}''_1$ defines a hyperplane $H_1(z)$ (intersection of $T\mathcal{F}''_1$ with the slices $\mathbb{D}^{2n}$ in $\mathbb{T}^2\times \mathbb{D}^2\times \mathbb{D}^{2n}$) in $\mathbb{R}^{2n}$ simply by making the $Z'$-component zero which according to the above is given by \[H_1(z)=\{A\in \mathbb{R}^{2n}: -\frac{1}{|z|^3}\Sigma(a_ix_i+b_iy_i)=(1-\psi)\}\] So if $Y\in H_1(z)$ then $-\frac{1}{|z|^3}\Sigma(x_i^2+y_i^2)=-\frac{1}{|z|}=(1-\psi)$ which is not possible as $(1-\psi)\geq 0$. So $Y\pitchfork H_1(z)$ and hence (as in lemma 2.1 of \cite{Bertelson}) $(\beta\oplus 0)\wedge \alpha+ d\alpha$ is nondegenerate on $(T\mathcal{F}''_1)_{(w,z)},\ for\ z\neq 0$ where $\beta$ is a one form on $\mathbb{T}^2\times \mathbb{D}^{q-2}$ such that $\beta(\psi Z')=1$ and $\beta=0$ on the orthogonal complement of $(\psi Z')$ in $T(\mathbb{T}^2\times \mathbb{D}^{q-2})$. So set $\omega''_1=\beta\wedge \alpha+d\alpha$. Observe that $\beta=0$ and outside the support of $\psi$, $\omega''_1=d\alpha=\omega'_1$. Moreover the pair $(\mathcal{F}'',\omega'')$ satisfies the properties of $(\mathcal{G},\gamma)$ of \ref{Main} only in the complement of the origin $\{z=0\}$.\\
 
 Now we shall change $\rho$ to $\tilde{\rho}$ and  $\alpha$ to $\tilde{\alpha}$ in a neighborhood of the origin $\{z=0\}$ such that by setting $\mathcal{G}=(id\times \tilde{\rho})^{-1}\tilde{\mathcal{F}}_1$ and $\gamma=(0\oplus \beta)\wedge \tilde{\alpha}+d\alpha$ we shall get the desired pair. So let us now define the $\tilde{\rho}$ and $\tilde{\alpha}$. \\
 
 First define $\bar{\rho}(z)=-(K/2)|z|^2+\Sigma_i(x_i+y_i)$, where $z=(x_1,y_1,...,x_n,y_n)$ and $K$ is a large positive number to be specified later. Now set \[\tilde{\rho}=B\bar{\rho}+(1-B)\rho\] where $B$ is a bump function which is $1$ on a $\delta/2$-neighborhood of $\{z=0\}$ and is equal to $0$ outside a $\delta$-neighborhood of $\{z=0\}$. Observe that \[d\bar{\rho}=(-Kx_i+1,-Ky_i+1)\] So we take $\delta<(1/K)$ and hence $\tilde{\rho}$ becomes regular. Now \[d\tilde{\rho}=(\bar{\rho}-\rho)dB+Bd\bar{\rho}+(1-B)d\rho\] So $d\tilde{\rho}(Y)$ remains negative on the complement of $\{z=0\}$ if we choose $K$ large. So now we need to define $\tilde{\alpha}$ in order to deal with the origin $\{z=0\}$. Observe that $\alpha=(1/2)\Sigma_i(x_idy_i-y_idx_i)$. We claim that it is enough to perturb the coefficients $x_i,y_i$'s of $dy_i,dx_i$'s in the expression of $\alpha$ suitably to get $\tilde{\alpha}$. We shall do the case $n=2$ i.e, on $\mathbb{D}^{q-2}\times \mathbb{T}^2\times \mathbb{D}^4$ the general case is same. Observe that on the $\delta/2$-neighborhood of the origin the hyperplane $H_1$ becomes \[H_1=\{(a_1,b_1,a_2,b_2):(1-Kx_1)a_1+(1-Ky_1)b_1+(1-Kx_2)a_2+(1-Ky_2)b_2=R(z)=1-\psi(z)\}\] and $Ker\alpha$ takes the form \[Ker\alpha=\{(a_1,b_1,a_2,b_2):x_1b_1-y_1a_1+x_2b_2-y_2a_2=0\}\] So we need to show that on $H_1\cap Ker\tilde{\alpha}$, $d\alpha=dx_1\wedge dy_1+dx_2\wedge dy_2$ is non-degenerate for a perturbation $\tilde{\alpha}$ of $\alpha$. In order to do so we compute $Ker\alpha\cap H_1$. This will suggest us the perturbation. So we need to solve the system of linear equation 
 \[
\begin{array}{rcl}
-y_1a_1+x_1b_1-y_2a_2+x_2b_2&=&0\\
(1-K\bar{x}_1)a_1+(1-K\bar{y}_1)b_1+(1-K\bar{x}_2)a_2+(1-K\bar{y}_2)b_2&=&R=(1-\psi(z))\\
\end{array} 
 \]
In the above the $\bar{x}_i,\bar{y}_i$'s are variables coming from $H_1$. As we shall perturb $\alpha$ we need to keep these different because we shall keep $H_1$ same but shall perturb $Y$ in order to perturb $\alpha$. Now solving we get $(a_1,b_1,a_2,b_2)$ where $b_2$ is given by \[(y_2(1-K\bar{y}_2)+x_2(1-K\bar{x}_2))^{-1}[y_2R-(y_2(1-K\bar{x}_1)-y_1(1-K\bar{x}_2))a_1-(y_2(1-K\bar{y}_1)+x_1(1-K\bar{x}_2))b_1]\] and $a_2$ is given by \[(y_2(1-K\bar{y}_2)+x_2(1-K\bar{x}_2))^{-1}[x_2R-(x_2(1-K\bar{x}_1)+y_1(1-K\bar{y}_2))a_1-(x_2(1-K\bar{y}_1)-x_1(1-K\bar{y}_2))b_1]\]

Now $d\alpha(a_1\partial_{x_1}+b_1\partial_{y_1}+a_2\partial_{x_2}+b_2\partial_{y_2},a'_1\partial_{x_1}+b'_1\partial_{y_1}+a'_2\partial_{x_2}+b'_2\partial_{y_2})=(a_1b'_1-a'_1b_1)+(a_2b'_2-a'_2b_2)$. Observe that if for $(a_1,b_1)\neq 0$, $(a_1b'_1-a'_1b_1)=0$ for all choice of $(a'_1,b'_1)$ then by choosing $b'_1=a_1$ and $a'_1=-b_1$ we get a contradiction. So this observation suggests that we need to compute $(a_2b'_2-a'_2b_2)$ and make the perturbation accordingly to make $(a_2b'_2-a'_2b_2)$ zero at the origin. $(a_2b'_2-a'_2b_2)$ is equal to \[(y_2(1-K\bar{y}_2)+x_2(1-K\bar{x}_2))^{-2}[R(a'_1-a_1)\{x_2y_2(1-K\bar{x}_1)-x_2y_1(1-K\bar{x}_2)+x_2y_2(1-K\bar{x}_1)\]\[-y_1y_2(1-K\bar{y}_2)\}+R(b'_1-b_1)\{x_2y_2(1-K\bar{y}_1)+x_1x_2(1-K\bar{x}_2)+x_2y_2(1-K\bar{y}_1)-x_1y_2(1-K\bar{y}_2)\}\]\[+(a_1b'_1-a'_1b_1)\{(x_2(1-K\bar{x}_1)-y_1(1-K\bar{y}_2))(y_2(1-K\bar{y}_1)+x_1(1-K\bar{x}_2))\]\[-(x_2(1-K\bar{y}_1)+x_1(1-K\bar{y}_2))(y_2(1-K\bar{x}_1)-y_1(1-K\bar{x}_2))\}]\]

Observe that $(a_2b'_2-a'_2b_2)=0$ for $\{x_1=0=y_1=x_2\}$ and $\{x_1=0=y_1=y_2\}$ but $x_2$ and $y_2$ can not be zero at the same time.\\

Now we know that $Y$ has the property $0\neq Y\pitchfork H_1$ (just put $x_i,y_i$'s in place of $a_i,b_i$'s in the equation defining $H_1$ and observe that $\delta<1/K$ and we are in a $\delta$-neighborhood of the origin) so it will continue to have this property if we $C^{\infty}$ perturb $Y$ slightly. Also observe that perturbing $Y$ would result a small perturbation in $\alpha$. So we do the perturbation now. We perturb $y_2$ by $f(y_2)$ so that $f(y_2)\neq 0$ on $(-\bar{\delta},\bar{\delta})$ for $0<\bar{\delta}<<\delta$. Hence according to the above observation $(a_2b'_2-a'_2b_2)=0$ at the origin and we are done.\\
 
 {\bf Acknowledgements}
 
I would like to thank the referee for pointing out some in-completions errors in an earlier version.


\begin{thebibliography}{99}
\bibitem{Bertelson}  Bertelson, M\'{e}lanie Foliations associated to regular Poisson structures. Commun. Contemp. Math. 3 (2001), no. 3, 441–456. (Reviewer: Edith Padr\'{o}n)
 
 \bibitem{Eliashberg}  Eliashberg, Y.; Mishachev, N. Introduction to the h-principle. Graduate Studies in Mathematics, 48. American Mathematical Society, Providence, RI, 2002. xviii+206 pp. ISBN: 0-8218-3227-1 (Reviewer: John B. Etnyre)
 
 \bibitem{Geiges}  Geiges, Hansj\"{o}rg An introduction to contact topology. Cambridge Studies in Advanced Mathematics, 109. Cambridge University Press, Cambridge, 2008. xvi+440 pp. ISBN: 978-0-521-86585-2 (Reviewer: John B. Etnyre)
 \bibitem{SchWeitzer}  Schweitzer, Paul A. Counterexamples to the Seifert conjecture and opening closed leaves of foliations. Ann. of Math. (2) 100 (1974), 386–400. (Reviewer: Robert Roussarie)
 \bibitem{Wilson}  Wilson, F. Wesley, Jr. On the minimal sets of non-singular vector fields. Ann. of Math. (2) 84 1966 529–536. (Reviewer: W. H. Gottschalk)
 
 \end{thebibliography}
\end{document}